\documentclass[amstex,12pt, amssymb]{article}

\usepackage{mathtext}
\usepackage[cp1251]{inputenc}
\usepackage[T2A]{fontenc}
\usepackage[dvips]{graphicx}
\usepackage{amsmath}
\usepackage{amssymb}
\usepackage{amsxtra}
\usepackage{latexsym}
\usepackage{ifthen}

\newcounter{lemma}[section]

\newcounter{corollary}[section]

\newcounter{remark}[section]

\newcounter{theorem}[section]

\newcounter{proposition}[section]

\newcounter{example}

\numberwithin{equation}{section}

\begin{document}

\markboth{\centerline{E.~SEVOST'YANOV, S.~SKVORTSOV, E.~PETROV}}
{\centerline{ON EQUICONTINUOUS FAMILIES ... }}

\def\cc{\setcounter{equation}{0}
\setcounter{figure}{0}\setcounter{table}{0}}

\overfullrule=0pt

%\normalsize\large

\author{{E.~SEVOST'YANOV, S.~SKVORTSOV, E.~PETROV}\\}

\title{
{\bf ON EQUICONTINUOUS FAMILIES OF MAPPINGS IN METRIC SPACES}}

\date{\today}
\maketitle

%\large
\begin{abstract} The article is devoted to the study of mappings with finite distortion in metric spaces.
Analogues of results relating to equicontinuity and normality of
families of quasiregular mappings are obtained. It is proved that
the indicated families are equicontinuous if the characteristic of
the mappings has a finite mean oscillation at each inner point, and
the maps omit a certain fixed continuum. An equicontinuity of
generalized quasiisometries on Riemannian manifolds is also
obtained.
\end{abstract}


\section{Introduction} There are well known conditions under which families of quasiregular mappings
are equicontinuous at a given point (see, for example,
\cite[Theorem~3.17]{MRV}, \cite[Theorem~2.4]{MSV} and
\cite[Corollary~III.2.7]{Ri}). Here we are talking about mappings
that omit some set of positive capacity, however, in the case of
homeomorphisms, it is possible to confine ourselves to a two-point
set (see~\cite[Theorem~19.2]{Va}). In particular, the following
result holds; see~\cite[Theorem~3.17]{MRV}.

\medskip
{\bf Theorem~(Martio-Rickman-V\"{a}is\"{a}l\"{a})}. {\sl Let
$n\geqslant 2,$ let $1\leqslant K<\infty$ and let $D$ be a domain in
${\Bbb R}^n.$ Then the family ${\cal Q}_{K, E}(D),$ consisting of
all $K$-quasiregular mappings $f:D\rightarrow \overline{{\Bbb
R}^n}\setminus E$ is equicontinuous (normal) in $D,$ if ${\rm
cap\,}E
>0.$}

\medskip
It is worth noting that the equicontinuity (normality) here should
be understood in terms of the so-called chordal metric,
see~\cite[Definition~12.1]{Va}. Without the requirement ${\rm
cap\,}E
>0,$ the result
is incorrect, as the example of the family of analytic functions
$f_n(z)=z^n,$ $z\in B(0, 2)\subset {\Bbb C},$ shows. Note that for
some classes of mappings with unbounded characteristic, similar
results were established in~\cite{Cr}, \cite{RS}, \cite{Sev$_1$} and
\cite{Sev$_2$}. In this publication, we show that in abstract metric
spaces there is a similar situation: families of mappings that omit
a certain continuum are equicontinuous at every point with
appropriate restrictions on the growth of their characteristics. We
separately consider Riemannian manifolds with isoperimetric
inequality, where this statement can be obtained under weaker
conditions. In the latter situation, the main result relates to
generalized quasiisometries, that is, mappings that distort the
$p$-modulus of families of paths according to the Poletsky
inequality, where $n-1<p<n, $ and $n$ is the dimension of the
manifold. As far as we know, the results obtained in the article are
new even in the quasiconformal case. This applies both to the case
of general metric spaces and to the situation of manifolds.

\medskip
Consider the following definition. Let $\overline{X}:=X\cup\infty,$
and let $h:\overline{X}\times\overline{X}\rightarrow {\Bbb R}$ be a
metric. We say that $h$ satisfies the {\it weak  sphericalization
condition,} if %$h$ is a metric on $\overline{X},$
$(\overline{X}, h)$ is a compact metric space and $h(x, y)\leqslant
d(x, y)$ for every $x, y\in X.$ A metric space $X$ is called a {\it
space admitting a weak sphericalization}, if there exists a metric
$h:\overline{X}\times\overline{X}\rightarrow {\Bbb R}$ satisfying
the weak sphericalization condition. There are a number of studies
directly devoted to the sphericalization of space, see, for example,
\cite{BB}, \cite{BHX}, \cite{LSh} and~\cite{DCX}.

\medskip
In what follows, $(X, d, \mu)$ and $\left(X^{\,\prime},
d^{\,\prime}, \mu^{\,\prime}\right)$ are metric spaces with metrics
$d$ and $d^{\,\prime}$ and locally finite Borel measures $\mu$ and
$\mu^{\,\prime},$ correspondingly. Let $G$ be a domain in $X,$ let
$Q:G\rightarrow [0, \infty]$ be a function measurable with respect
to measure $\mu$ and $p,q \geqslant 1$ be fixed numbers. Put
$$A(\zeta_0, r_1, r_2)=\left\{x\in X: r_1<d(x, \zeta_0)<r_2\right\}\,,$$
$$S_i=S(\zeta_0, r_i)=\{x\in X: d(x, \zeta_0)=r_i\}\,.$$
Similarly to \cite[Ch.~7]{MRSY}, a mapping $f:G\rightarrow
X^{\,\prime}$ is called a {\it ring $Q$-mapping at a point $x_0\in
\overline{G}$ with respect to $(p, q)$-moduli}, if for some
$r_0=r_0(x_0)>0$ and for all $0<r_1< r_2<r_0$ the inequality
\begin{equation}\label{eq1}
M_p(f(\Gamma(S(\zeta_0, r_1), S(\zeta_0, r_2), A)))\leqslant
\int\limits_A Q(x)\cdot \eta^q(d(x, \zeta_0))\,d\mu(x)
\end{equation}
holds for any measurable function
$\eta:(r_1, r_2)\rightarrow [0, \infty]$ with
$\int\limits_{r_1}^{r_2}\eta(r)\,dr\geqslant 1.$
If $p=q,$ then for brevity we will call such mappings {\it ring
$Q$-mappings with respect to a $p$-modulus}. For many well-known
mapping classes, estimates of the form~(\ref{eq1}) are established.
In particular, the class of ring $Q$-mappings includes analytic
functions with $X=X^{\,\prime}={\Bbb C},$ $Q(x)\equiv 1,$ $p=q=2,$
and also quasiregular (quasiconformal) mappings with
$X=X^{\,\prime}={\Bbb R}^n,$ $Q(x)\leqslant K=const$ and $p=q=n.$
There are many mappings with unbounded characteristic for which the
estimates~(\ref{eq1}) are also satisfied. In particular,
homeomorphisms $f\in W_{\rm loc}^{1, n},$ for which $f^{\,-1}\in
W_{\rm loc}^{1, n};$ in this case, $Q=K_I(x, f)$ and $p=q=n,$ where
$K_I(x, f)$ is the inner dilatation of the map $f$ at the point $x$
(see~\cite[Theorems~8.1, 8.6]{MRSY}).

\medskip
Similarly to~\cite[Section~13.4]{MRSY}, we say that a function
$\varphi:G\rightarrow{\Bbb R}$ has {\it finite mean oscillation at a
point $x_{0}\in\overline{G}$}, abbreviated $\varphi \in FMO(x_{0})$,
if
\begin{equation}\label{eq13.4.111} \overline{\lim\limits_{\varepsilon\rightarrow
0}}\,\, \,\frac{1}{\mu(B(x_{0},\varepsilon))}
\int\limits_{B(x_{0},\varepsilon)}|\varphi(x)-\overline{\varphi}_{\varepsilon}|\,\,d\mu(x)<\infty\,,
\end{equation}
where $\overline{\varphi}_{\varepsilon}
=\frac{1}{\mu(B(x_{0},\varepsilon))}
\int\limits_{B(x_{0},\varepsilon)}\varphi(x)\,\,d\mu(x)$ is the mean
value of the function $\varphi(x)$ over the ball
$B(x_{0},\varepsilon)=\{x\in G: d(x,x_0)<\varepsilon\}$ with respect
to the measure $\mu$. Here the condition (\ref{eq13.4.111}) includes
the assumption that $\varphi$ is integrable with respect to the
measure $\mu$ over the ball $B(x_0,\varepsilon)$ for some
$\varepsilon>0.$ Metric measure spaces where the inequalities
$\frac{1}{C}R^{n}\leqslant \mu(B(x_0, R))\leqslant CR^{n}$
hold for a constant $C\geqslant 1$, every $x_0\in X$  and all
$R<{\rm diam}\,X$, are called {\it Ahlfors $n$-regular.} Following
\cite[section~7.22]{He}, given a real-valued function $u$ in a
metric space $X,$ a Borel function $\rho\colon X\rightarrow [0,
\infty]$ is said to be an {\it upper gradient} of a function
$u:X\rightarrow{\Bbb R}$ if $$|u(x)-u(y)|\leqslant
\int\limits_{\gamma}\rho\,|dx|$$ for each rectifiable path $\gamma$
joining $x$ and $y$ in $X.$ Let $(X, \mu)$ be a metric measure space
and let $1\leqslant p<\infty.$
We say that $X$ admits {\it a $(1; p)$-Poincar\'{e} inequality} if
there is a constant $C\geqslant 1$ and $\tau\geqslant 1$ such that
$$\frac{1}{\mu(B)}\int\limits_{B}|u-u_B|\,d\mu(x)\leqslant C\cdot({\rm
diam\,}B)\left(\frac{1}{\mu(\tau B)} \int\limits_{\tau B}\rho^{\,p}
\,d\mu(x)\right)^{1/p}$$
for all balls $B$ in $X,$ for all bounded continuous functions $u$
on $B,$ and for all upper gradients $\rho$ of $u.$ Let $X$ and $Y$
be metric spaces. A mapping $f:X\rightarrow Y$ is {\it discrete} if
$f^{\,-1}(y)$ is discrete for all $y\in Y$ and $f$ is {\it open} if
$f$ maps open sets onto open sets.

\medskip
Following~\cite{SM}, we agree in the following terminology. Given
$2\leqslant\alpha<\infty$ and $1\leqslant q\leqslant \alpha,$ the
space $X=(X, d, \mu)$ is called {\it $(\alpha, q)$-admissible
source}, if $(X, d, \mu)$ be locally compact and locally path
connected Ahlfors $\alpha$-regular metric space. Similarly, given
$2\leqslant\alpha^{\,\prime}<\infty$ and
$\alpha^{\,\prime}-1<p\leqslant\alpha^{\,\prime},$ the space
$X^{\,\prime}=(X^{\,\prime}, d^{\,\prime}, \mu^{\,\prime})$ is
called {\it $(\alpha^{\,\prime}, p)$-admissible target}, if
$(X^{\,\prime}, d^{\,\prime}, \mu^{\,\prime})$ be an Ahlfors
$\alpha^{\,\prime}$-regular, path connected, locally connected and
locally compact metric space which supports $(1; p)$-Poincar\'{e}
inequality.

\medskip
We introduce the class of mappings to which our studies are devoted.
Let $K$ be some continuum in $X^{\,\prime}.$ Denote by $\frak{F}^{p,
q}_{\zeta_0, Q}(G, X^{\,\prime}, K)$ the family of all open discrete
ring $Q$-mappings $f:G\rightarrow X^{\,\prime}\setminus K$ at the
point $\zeta_0\in G$ with respect to the $(p, q)$-moduli. The
following statement generalizes Montel's theorem on the normality of
families of flat analytic functions, see~\cite[Theorems~19.4 and
20.5]{Va} and~\cite[Theorem~5.11]{RS}.

\medskip
\begin{theorem}\label{th1}{\sl\,
Let $2\leqslant\alpha, \alpha^{\,\prime}<\infty,$ let
$\alpha^{\,\prime}-1<p\leqslant\alpha^{\,\prime}$ and $1\leqslant
q\leqslant \alpha,$ let $(X, d, \mu)$ be an $(\alpha, q)$-admissible
source and let $(X^{\,\prime}, d^{\,\prime}, \mu^{\,\prime})$ be an
$(\alpha^{\,\prime}, p)$-admissible target. Assume that $Q\in
FMO(\zeta_0),$ and that $X^{\,\prime}$ admits a weak
sphericalization. Then the family $\frak{F}^{p, q}_{\zeta_0, Q}(G,
X^{\,\prime}, K)$ is equicontinuous at $\zeta_0$ in the sense of the
space $(\overline{X^{\,\prime}}, h).$
}
 \end{theorem}

\medskip
A metric space is said to be {\it proper} if its closed balls are
compact sets. A space $X$ is called {\it Ptolemaic} if, for every
$x, y, z, t\in X$, we have that
%
%\begin{equation}\label{eq1B}
$$d(x, z)d(y, t) + d(x, t)d(y, z) - d(x, y)d(z, t) \geqslant 0\,.$$
%\end{equation}
%
It can be shown that proper Ptolemaic spaces admit weak
sphericalization, where as $ h $ it is necessary to consider the
metric
\begin{equation}\label{eq2.2.1}
h_{x_0}(x,y):=\frac{d(x, y)}{\sqrt{1+d^2(x, x_0)}\sqrt{1+d^2(y,
x_0)}}\,,
\end{equation}
$$h_{x_0}(x, \infty)=\frac{1}{\sqrt{1+d^2(x_0, x)}}\,,$$
and $x_0\in X$ is some fixed point (see~\cite[Lemma~5.4]{SM}). If
$X={\Bbb R}^n$ and $x_0=0,$ then the metric $h_{x_0},$ defined
in~(\ref{eq2.2.1}), determines the distance between the projections
of the points on the Riemann sphere (in this case, it is called a
{\it chordal} metric, see~\cite[Definition~12.1]{Va}). Note that the
space ${\Bbb R}^n$ is Ptolemaic, see~\cite[Proposition~10.9.2]{B},
cf.~\cite{DP}. We define~$\frak{F}^{p, q}_{Q}(G, X^{\,\prime}, K)$
as a family of mappings $f:G\rightarrow X^{\,\prime},$ belonging to
$\frak{F}^{p, q}_{\zeta_0, Q}(G, X^{\,\prime}, K)$ for each
$\zeta_0\in G.$ Taking into account the above, from
Theorem~\ref{th1} and the Arzel–Ascoli theorem
(see~\cite[Theorem~20.4]{Va}), we obtain the following statement.

\medskip
\begin{corollary}\label{cor1}
{\sl\, Let $2\leqslant\alpha, \alpha^{\,\prime}<\infty,$ let
$\alpha^{\,\prime}-1<p\leqslant\alpha^{\,\prime}$ and $1\leqslant
q\leqslant \alpha,$ let $(X, d, \mu)$ be an $(\alpha, q)$-admissible
source and let $(X^{\,\prime}, d^{\,\prime}, \mu^{\,\prime})$ be an
$(\alpha^{\,\prime}, p)$-admissible target. Assume that $Q\in
FMO(\zeta_0)$ for all $\zeta_0\in G,$ $(X, d, \mu)$ is separable and
$(X^{\,\prime}, d^{\,\prime}, \mu^{\,\prime})$ is Ptolemaic and
proper. Then $\frak{F}^{p, q}_{Q}(G, X^{\,\prime}, K)$ is a normal
family of mappings in $G$ in the sense of the space
$(\overline{X^{\,\prime}}, h_{x_0}).$ }
\end{corollary}

\medskip
\section{The main lemma and proof of
Theorem~\ref{th1}}\setcounter{equation}{0}

As usually, given a path $\gamma:[a, b]\rightarrow X,$ we set
$$|\gamma|:=\{x\in X: \exists\,t\in[a, b]:
\gamma(t)=x\}\,.$$
Recall that a pair $E=\left(A,\,C\right),$ where $A$ is an open set
in $X,$ and $C\subset A$ is a compact set, is called {\it condenser}
in $X.$ Given $p\geqslant 1,$ a quantity
\begin{equation}\label{eq28}
{\rm cap}_p\,E=M_p(\Gamma_E)
\end{equation}
is called {\it $p$-capacity of $E,$} where $\Gamma_E$ be the family
of all paths of the form $\gamma\colon [a,\,b)\rightarrow A$ with
$\gamma(a)\in C$ and $|\gamma|\cap\left(A\setminus
F\right)\ne\varnothing$ for every compact set $F\subset A.$ The
following most important statement is proved in
~\cite[Lemma~3.3]{SM}.

\medskip
\begin{lemma}\label{lem4}
{\sl\,Let $\alpha\geqslant 2,$ let $\alpha-1<p\leqslant\alpha,$ and
let $X$ be $(\alpha, p)$-admissible target. Suppose that $X$ admits
a weak sphericalization and $F$ is nondegenerate continuum in $X.$
Now, for every $a>0$ there exists $\delta>0$ for which the condition
\begin{equation*}\label{eq5}{\rm cap}_p\,\left(X\setminus F,\, C\right)\geqslant \delta\end{equation*}
holds for any continuum $C\subset X\setminus F$ with $h(C)\geqslant
a.$
 }
\end{lemma}

The following lemma contains the statement of Theorem~\ref{th1} in
the case when the function $Q$ from the definition of the class
$\frak{F}^{p, q}_{\zeta_0, Q}(G, X^{\,\prime}, K)$ satisfies the
most general formal conditions.

\medskip
\begin{lemma}\label{lem1}{\sl\,
Let $2\leqslant\alpha, \alpha^{\,\prime}<\infty,$ let
$\alpha^{\,\prime}-1<p\leqslant\alpha^{\,\prime}$ and $1\leqslant
q\leqslant \alpha,$ let $(X, d, \mu)$ be an $(\alpha, q)$-admissible
source and let $(X^{\,\prime}, d^{\,\prime}, \mu^{\,\prime})$ be an
$(\alpha^{\,\prime}, p)$-admissible target. Suppose that
$\varepsilon_0>0$ and $\psi\colon(0, \varepsilon_0)\rightarrow
[0,\infty]$ is a Lebesgue measurable function that satisfies the
following condition: for every $\varepsilon_2\in(0, \varepsilon_0]$
there is $\varepsilon_1\in (0, \varepsilon_2]$ such that for every
$\varepsilon\in (0,\varepsilon_1)$
\begin{equation}\label{eq3AB} 0<I(\varepsilon, \varepsilon_0):=
\int\limits_{\varepsilon}^{\varepsilon_0}\psi(t)\,dt <
\infty\qquad\forall\,\,\varepsilon\in(0, \varepsilon_1)\,.
\end{equation}
Suppose also that
\begin{equation} \label{eq3.7C}
\int\limits_{\varepsilon<d(x, \zeta_0)<\varepsilon_0}
Q(x)\cdot\psi^q(d(x, \zeta_0))\,d\mu(x)=o(I^{\,q}(\varepsilon,
\varepsilon_0))\,,
\end{equation}
and that $X^{\,\prime}$ admits a weak sphericalization. Then the
family $\frak{F}^{p, q}_{\zeta_0, Q}(G, X^{\,\prime}, K)$ is
equicontinuous at $\zeta_0$ in the sense of the space
$(\overline{X^{\,\prime}}, h).$ }
\end{lemma}

\medskip
\begin{proof}
Since the space $X$ is locally compact, we may assume that
$\overline{B(\zeta_0, \varepsilon_0)}$ is a compact in $G.$ Let
$\beta: [a,\,b)\rightarrow X^{\,\prime}$ and let
$x\in\,f^{\,-1}(\beta(a)).$ Since $X$ and $X^{\,\prime}$ are locally
compact, $X$ is locally connected, and $f:D \rightarrow
X^{\,\prime}$ is discrete and open, there is a maximal $f$-lifting
of $\beta$ starting at $x,$ see~\cite[Lemma~2.1]{SM}. In this case,
, setting $S_1=S(\zeta_0, \varepsilon),$ $S_2=S(\zeta_0,
\varepsilon_0),$ and arguing similarly to the proof
of~\cite[Lemma~3]{Sev$_3$}, we obtain that
\begin{equation}\label{eq3C}
M_p(\Gamma(f(\overline{B(\zeta_0, \varepsilon)}),
\partial f(B(\zeta_0, \varepsilon_0)), X^{\,\prime}))=\alpha(\varepsilon)\rightarrow 0\,,\qquad\varepsilon\rightarrow
0\,,
\end{equation}
$0<\varepsilon<\varepsilon^{\,\prime}_0,$ where $A=A(\zeta_0,
\varepsilon, \varepsilon_0)=\left\{x\in X: \varepsilon<d(x,
\zeta_0)<\varepsilon_0\right\},$ and $I$ is defined
in~(\ref{eq3AB}). Consider a condenser ${\cal E}=(A, C),$
$A=B\left(\zeta_0, \varepsilon_0\right),$ $C=\overline{B(\zeta_0,
\varepsilon)}.$ Then, in terms of the capacity of the condenser
${\cal E},$ the relation~(\ref{eq3C}) can be rewritten as
\begin{equation}\label{eq3CA}
{\rm cap}_p\,f{\cal E}=\alpha(\varepsilon)\rightarrow 0\,,
\qquad\varepsilon\rightarrow 0\,.
\end{equation}
On the other hand, by Lemma~\ref{lem4} for each  $a>0$ there exists
$\delta=\delta(a)$ such that the relation
\begin{equation}\label{eq5A}
{\rm cap}_p\,\left(X^{\,\prime}\setminus K,\, C\right)\geqslant
\delta\,.\end{equation}
is fulfilled for any arbitrary continuum $C\subset
X^{\,\prime}\setminus K,$ satisfying the condition $h(C)\geqslant
a.$ The estimate~(\ref{eq3CA}) implies that for $\delta=\delta(a)$
there is $\varepsilon_*=\varepsilon_*(a)$ such that
\begin{equation}\label{eq18}
{\rm cap}_p\,f{\cal E}\leqslant \delta\quad \forall\,\,
\varepsilon\,\in \left(0,\,{\varepsilon_*}(a)\right)\,.
\end{equation}
Using the relation~(\ref{eq18}), we obtain that
$${\rm cap}_p\,\left(X^{\,\prime}\setminus
K,\,f\left(\overline{B(\zeta_0,\,\varepsilon)}\right)\right)\,\leqslant$$
\begin{equation*}\label{eq19}
\leqslant {\rm
cap}_p\,\left(f\left({B(\zeta_0,\,\varepsilon_0)}\right)\,,\,
f\left(\overline{B(\zeta_0,\,\varepsilon)}\right)\right)\,<\,\delta
\end{equation*}
for $\varepsilon\,\in \left(0,\,{\varepsilon_*}(a)\right).$ Then it
follows from~(\ref{eq5A}) that
$h\left(f\left(\overline{B(\zeta_0,\,\varepsilon)}\right)\right)<a.$
Finally, for any $a>0$ there exists $\varepsilon_*=\varepsilon_*(a)$
such that
$h\left(f\left(\overline{B(\zeta_0,\,\varepsilon)}\right)\right)<a$
whenever $\varepsilon\in \left(0,\,{\varepsilon_*}(a)\right).$ The
lemma is proved.~$\Box$
\end{proof}

\medskip
{\it The proof of Theorem~\ref{th1}} immediately follows from
Lemma~\ref{lem1} and \cite[Proposition~2]{Sev$_3$}.~$\Box$

\section{Auxiliary information from the theory of metric spaces and manifolds}
\setcounter{equation}{0}

As we said above, another important result of the article relates to
Riemannian manifolds. Recall that the length of a piecewise smooth
path $\gamma=\gamma(t),$ $t\in [a, b],$ joining points connecting
points $\gamma(a)=M_1\in {\Bbb M}^n$ and $\gamma(b)=M_2\in {\Bbb
M}^n$ on the Riemannian manifold ${\Bbb M}^n$ is determined by the
relation
\begin{equation}\label{eq2A}
l(\gamma):=\int\limits_{a}^{b}\sqrt{\sum\limits_{i,j=1}^ng_{ij}(x(t))\frac{dx^i}{dt}\frac{dx^j}{dt}}\,dt\,,
\end{equation}
where $g=g_{ij}(x)$ is a smooth positive definite tensor of type
$(0, 2) $ on a manifold (Riemannian metric). In other words,
$g=g_{ij}(x)$ is a system of matrices that in various coordinate
systems are connected by the relation
$'g_{ij}(x)=\sum\limits_{k,l=1}^n g_{kl}(y(x))\frac{\partial
y^k}{\partial x^i} \frac{\partial y^l}{\partial x^j}.$ {\it Geodetic
distance} between points $p_1$ and $p_2\in{\Bbb M}^n$ will be called
the shortest length of all piecewise smooth paths in ${\Bbb M}^n,$
joining the points $p_1$ and $p_2.$ Geodesic distance between points
$p_1$ and $p_2$ will be denoted by $d(p_1, p_2).$ In these terms, it
also makes sense
\begin{equation}\label{eq6}
l(\gamma):=\sup\sum\limits_{i=1}^m d(\gamma(t_i),
\gamma(t_{i-1}))\,,
\end{equation}
where $\sup$ is taken over all possible partitions $a=t_0\leqslant
t_1\leqslant\ldots\leqslant t_m:=b.$ So, we have ''two definitions''
for $l(\gamma),$ (\ref{eq2A}) and (\ref{eq6}), one of which,
corresponding to~(\ref{eq6}), suits for arbitrary curves, not just
piecewise smooth. In this regard, the following remark must be made.

\medskip
\begin{remark}\label{rem1}
For piecewise smooth paths $\gamma:[a, b]\rightarrow {\Bbb M}^n,$
belonging to the connected Riemannian manifold ${\Bbb M}^n,$ the
quantities defined by the relations (\ref{eq2A}) and (\ref{eq6}),
coincide (see~\cite[Theorem~2.2]{Bur}).
\end{remark}

\medskip
\begin{lemma}\label{lem5}
{\sl Let $(X, d)$ be a proper metric space. Then $(X, d)$ is
complete.}
\end{lemma}

\medskip
\begin{proof}
Indeed, let $x_m,$ $m=1,2,\ldots ,$ be a fundamental sequence in
$X.$ Then for $\varepsilon=1$ there is a positive integer $M>1$ such
that $d(x_m, x_l)<1$ for all $m, l\geqslant M.$ In this case,
$x_l\in B(x_{M}, 1),$ $l\geqslant M.$ Since the space $X$ is proper,
the closure of any ball $B(a, r)$ is a compact in $X.$ Moreover,
since all elements of the sequence $\{x_m\}_{m=1}^{\infty}$ (except
a finite number) belong to the ball $B(x_M, 1),$ there exists a
subsequence $x_{m_k}$ of the sequence $x_m$ converging to some
element $x_0\in X$ as $k\rightarrow\infty.$

Since the sequence $x_m,$ $m=1,2,\ldots$ is fundamental, for an
arbitrary  $\varepsilon>0$ there is a number $N=N(\varepsilon)$ such
that $d(x_m, x_l)<\varepsilon/2$ for all $m,l> N.$ Moreover, due to
the convergence of $x_{m_k}$ to $x_0$ as $k\rightarrow\infty,$ for
some $K=K(\varepsilon)$ and all $k> K,$ we obtain that
$d(x_{m_k},x_0)<\varepsilon/2.$ Using these two inequalities, by the
triangle inequality we obtain that $d(x_l, x_0)\leqslant d(x_l,
x_{m_k})+d(x_{m_k}, x_0)<\varepsilon/2+\varepsilon/2=\varepsilon$
for $l, k>\max\{K, N \},$ i.e., the sequence $x_m$ converges as
$m\rightarrow\infty$ to $x_0.$~$\Box$
\end{proof}

\medskip
Recall that a metric space $(X, d)$ is called {\it geodesic} if any
two points $x_1, x_2\in X$ can be joined by a rectifiable path
$\gamma:[0, 1]\rightarrow X,$ $\gamma(0)=x_1$ and $\gamma(1)=x_2,$
whose length coincides with $d(x_1, x_2).$ The following lemma
holds, cf.~\cite[Remark~9.11]{He}.

\medskip
\begin{lemma}\label{lem3}
{\sl Let ${\Bbb M}^n,$ $n\geqslant 2,$ be a Riemannian manifold with
a geodesic metric $d,$ which is a connected and proper space. Then
the space ${\Bbb M}^n$ is a geodesic.}
\end{lemma}

\medskip
\begin{proof}
By Lemma~\ref{lem5}, ${\Bbb M}^n$ is a complete Riemannian manifold,
which means that for any two points $p_1, p_2 \in {\Bbb M}^n$ there
is a piecewise smooth path $\gamma$ path joining these points whose
length coincides with $d(p_1, p_2)$ in the sense of (\ref{eq2A}),
see~\cite[Corollary~6.15]{Lee}. Moreover, this length coincides with
the length of $\gamma,$ understood in the sense of a metric space
and defined by the formula~(\ref{eq6}) (see Remark~\ref{rem1}). The
lemma is completely proved.~$\Box$
\end{proof}

\medskip
The following result holds (see \cite[Proposition~4.7]{AS}).

\medskip
\begin{proposition}\label{pr2}
{\sl Let $X$ be a $Q$-Ahlfors regular metric measure space that
supports $(1; p)$-Poincar\'{e} inequality for some $p>1$ such that
$Q-1<p\leqslant Q.$ Then there exists a constant $C>0$ having the
property that, for $x\in X,$ $R>0$ and continua $E$ and $F$ in $B(x,
R),$
$$M_p(\Gamma(E, F, X))\geqslant \frac{1}{C}\cdot\frac{\min\{{\rm diam}\,E, {\rm diam}\,F\}}{R^{1+p-Q}}\,.$$}
\end{proposition}

\section{On normal families of generalized quasiisometries in spaces
with isoperimetric inequality} \setcounter{equation}{0}

We consider separately the possibility of $n-1<p <n,$ that is, we
exclude the case $p=n$ from consideration. If we are talking about
$n$-dimensional Euclidean space, the corresponding families of
mappings are equicontinuous even without the assumption of omitting
some set $E,$ see~\cite[Theorem~1.1]{GSS}. We now show that
something similar is true for Riemannian manifolds with
isoperimetric inequality. These inequalities are well known in the
Euclidean case (see, for example,~\cite[Section~3.2.43]{Fe}). Here
we restrict ourselves to the consideration of families of
homeomorphisms. Everywhere below $d$ and $d_*$ are geodesic
distances on connected Riemannian manifolds ${\Bbb M}^n$ and ${\Bbb
M}^n_*,$ respectively, and $dv$ and $dv_*$ are volume elements on
these manifolds. The concepts associated with their study, we
consider known (they can be found, for example, in~\cite{Lee},
cf.~\cite{IS$_1$}).

\medskip
In what follows, $p$-capacity of the condenser $E=(A, C)$ in ${\Bbb
M}^n$ (or ${\Bbb M}^n_*,$ depending on the context) is defined by
the relation~(\ref{eq28}). It is worth noting that on manifolds this
definition coincides with the well-known definition of capacity
through the integral of the gradient (see, for example,
\cite[Proposition~II.10.2 and Remark~II.10.8]{Ri}). We now require
the following additional condition on the Riemannian manifold ${\Bbb
M}^n_*$ with the geodesic distance $d_*$ and the measure of volume
$v_*.$ Denoting by $\mathcal{H}^{n-1}(B)$ the $(n-1)$-dimensional
Hausdorff measure of the measurable set $B\subset {\Bbb M}^n_*,$ we
assume that for any open set $A\subset {\Bbb M}^n_*$ with compact
closure and smooth boundary, the condition
\begin{equation}\label{eq3}
\mathcal{H}^{n-1}(\partial A)\geqslant \frac{c}{(v_*(A))^{1/n-1}}
\end{equation}
holds for some constant $c>0.$ The condition~(\ref{eq3}) is called
{\it isoperimetric inequality} on ${\Bbb M}^n_*.$ In this case,
by~\cite[inequalities~(1.7) and~(4.1)]{Gr}, the estimate
\begin{equation}\label{eq3A}
{\rm cap}_p (A, C)\geqslant C^{\,\prime}\cdot (v_*(C))^{1-p/n}\,.
\end{equation}
holds for any condenser $E=(A, C)$ in ${\Bbb M}^n_*.$

Let $G$ be a domain of a Riemannian manifold ${\Bbb M}^n,$
$n\geqslant 2,$ and let $Q:G\rightarrow [0, \infty]$ be a measurable
function. Denote by $\frak{F}^p_{\zeta_0, Q}(G, {\Bbb M}^n_*)$ the
family of all ring $Q$-homomorphisms $f:G\rightarrow X^{\,\prime}$
at $\zeta_0\in G$ with respect to the $p$-modulus. The following
result holds, cf.~\cite[Lemma~2.4]{GSS}.

\medskip
\begin{lemma}\label{lem2}{\sl\,
Let $Q:G\rightarrow (0, \infty]$ be a function measurable with
respect to measure $v,$ $n-1<p<n,$ and the Riemannian manifold
${\Bbb M}_*^n$ be considered as a metric space with a geodesic
metric $d_*,$ is $(n, p)$-admissible target. In addition, assume
that ${\Bbb M}_*^n$ is a proper space in which the isoperimetric
inequality~(\ref{eq3}) is fulfilled. Suppose that the
conditions~(\ref{eq3AB})--(\ref{eq3.7C}) are satisfied, where $q=p.$
Then the family $\frak{F}^p_{\zeta_0, Q}(G, {\Bbb M}^n_*)$ is
equicontinuous at $\zeta_0,$ where equicontinuity is understood in
the sense of the space $({\Bbb M}^n_*, d_*).$
}
 \end{lemma}

\medskip
\begin{proof}
\textbf{I.} Suppose the contrary, namely, suppose that the family
$\frak{F}^p_{\zeta_0, Q}(G, {\Bbb M}^n_*)$ is not equicontinuous at
some point $\zeta_0\in G.$ Then there is $\delta_0>0$ and sequences
$\zeta_m\in G,$ $f_m\in \frak{F}^p_{\zeta_0, Q}(G, {\Bbb M}^n_*)$
such that $\zeta_m\rightarrow \zeta_0$ as $m\rightarrow\infty$ and
\begin{equation}\label{eq4A}
d_*(f_m(\zeta_m), f_m(\zeta_0))\geqslant\delta_0\,.
\end{equation}
In this case, we fix two points $x_1, x_2\in{\Bbb M}^n_*,$ $x_1\ne
x_2$ in an arbitrary way. Without loss of generality, we may assume
that $B_1:=B(x_1, \delta_0/3)$ and $B_2:=B(x_2, \delta_0/3)$ are
compact in ${\Bbb M}^n_*$ and $\delta_0<d_*(x_1, x_2).$ (Otherwise,
we may put $\delta_0^{\,*}:=d_*(x_1, x_2)/2.$ Then we obtain an
inequality similar to~(\ref{eq4A}), in which $\delta_0$ is replaced
by $\delta_0^{\,*}.$ In this situation, $\delta_0^{\,*}<d_*(x_1,
x_2)$). By the triangle inequality
\begin{equation}\label{eq11}
d_*(\overline{B_1}, \overline{B_2})\geqslant \delta_0/3\,.
\end{equation}
Since ${\Bbb M}^n_*$ is Ahlfors $n$-regular, there exists some
constant $\widetilde{C}>0,$ depending only on ${\Bbb M}^n_*$ such
that
\begin{equation}\label{eq12}
v_*(\overline{B_i})\geqslant \widetilde{C}\cdot \delta_0^n,\quad
i=1,2\,.
\end{equation}

\medskip
\textbf{II.} Consider the condenser ${\mathcal E}=(A,\,C),$ where
$A=B\left(\zeta_0,\,\varepsilon_0\right),$
$C=\overline{B(\zeta_0,\,\varepsilon)}$ and
$0<\varepsilon<\varepsilon_0.$ Without loss of generality, we may
assume that $C$ belongs to some normal neighborhood of $\zeta_0$
and, in particular, the sets $B(\zeta_0, \varepsilon)=\{x\in{\Bbb
M}_*^n: d(x,\zeta_0)<\varepsilon\}$ and $S(\zeta_0,
\varepsilon)=\{x\in{\Bbb M}_*^n: d(x, \zeta_0)=\varepsilon \}$ are
path connected for all $\varepsilon\in(0, \varepsilon_0).$ (The
definition of normal neighborhoods of a point can be found, for
example, in~\cite[Proposition~5.11]{Lee}).

Arguing as in the proof of Lemma~\ref{lem1}, we conclude that the
relation~(\ref{eq3C}) holds for $p=q$ and $X^{\,\prime}={\Bbb
M}_*^n.$ This implies that
\begin{equation*}\label{eq16}
{\rm cap}_p\,f({\mathcal E})\leqslant
\alpha(\varepsilon)\quad\forall\,\,\varepsilon\,\in
(0,\,\varepsilon_0),\quad f\in\frak{F}^p_{\zeta_0, Q}(G, {\Bbb
M}^n_*),\quad
\alpha(\varepsilon)\,\stackrel{\varepsilon\,\rightarrow
0}{\rightarrow}\,0\,.
\end{equation*}
On the other hand, by~(\ref{eq3}), in view of the above remarks,
(\ref{eq3A}) also holds. Consequently,
\begin{equation*}
\alpha(\varepsilon)\geqslant{\rm cap}_p\,f({\mathcal E})\geqslant
C^{\,\prime}\cdot\left[v_*(f(C))\right]^{\frac{n-p}{n}}\quad
\forall\,\, f\in\frak{F}^p_{\zeta_0, Q}(G, {\Bbb M}^n_*) \,.
\end{equation*}
From the last relation it follows that
\begin{equation*}\label{eq6A}
v_*(f(C))\leqslant \alpha_1(\varepsilon) \quad \forall\,\,
f\in\frak{F}^p_{\zeta_0, Q}(G, {\Bbb M}^n_*)\,,
\end{equation*}
where $\alpha_1$ is some function such that
$\alpha_1(\varepsilon)\rightarrow 0$ as $\varepsilon\rightarrow 0.$
Therefore,
\begin{equation}\label{eq10}
v_*(f(\overline{B(\zeta_0, \varepsilon_1)}))\leqslant
\widetilde{C}\cdot \delta_0^n/2 \quad \forall\,\,
f\in\frak{F}^p_{\zeta_0, Q}(G, {\Bbb M}^n_*)
\end{equation}
for some $\varepsilon_1>0.$

\medskip
\textbf{III.} Denote $A_1:=B(\zeta_0, \varepsilon_1),$
$C_1:=\overline{B(\zeta_0, \varepsilon)}$ and consider another
condenser  ${\mathcal E}_1=(A_1, C_1)=(B(\zeta_0, \varepsilon_1),
\overline{B(\zeta_0, \varepsilon)}),$ where
$0<\varepsilon<\varepsilon_1.$ From the
conditions~(\ref{eq3AB})--(\ref{eq3.7C}), taking into account the
analogue of the Fubini theorem on Riemannian manifolds, we conclude
that the integral in~(\ref{eq3.7C}) does not tend to zero for small
values of $\varepsilon>0,$ see~\cite[Theorem~2.1]{IS$_1$}. At the
same time, we also took into account that $Q(x)>0$ almost everywhere
by the condition of the lemma. Therefore, $I(\varepsilon,
\varepsilon_0)\rightarrow \infty$ as $\varepsilon\rightarrow 0.$ It
follows that the relations~(\ref{eq3AB})--(\ref{eq3.7C}) remain
valid if we replace $\varepsilon_0$ with $\varepsilon_1\in(0,
\varepsilon_0).$

Then, arguing as in paragraph \textbf{II}, we conclude that
$${\rm cap}_p\,f({\mathcal E}_1)={\rm cap}_p\,(f(B(\zeta_0,
\varepsilon_1)), f(\overline{B(\zeta_0, \varepsilon)}))\leqslant$$
\begin{equation}\label{eq16A}
\leqslant \alpha_2(\varepsilon)\quad\forall\,\,\varepsilon\,\in
(0,\,\varepsilon_0),\quad\forall\,\, f\in\frak{F}^p_{\zeta_0, Q}(G,
{\Bbb M}^n_*),\quad
\alpha_2(\varepsilon)\,\stackrel{\varepsilon\,\rightarrow
0}{\rightarrow}\,0\,.
\end{equation}

\medskip
\textbf{IV.} Note that both in the ball $\overline{B_1},$ and in the
ball $\overline{B_2},$ there are points $z_m\in \overline{B_1},$
$w_m\in \overline{B_2},$ such that $z_m, w_m\not\in f_m(A_1),$ which
directly follows from the relations (\ref{eq12}) and (\ref{eq10}),
see Figure~\ref{fig1}.
\begin{figure}
  \centering\includegraphics[width=350pt]{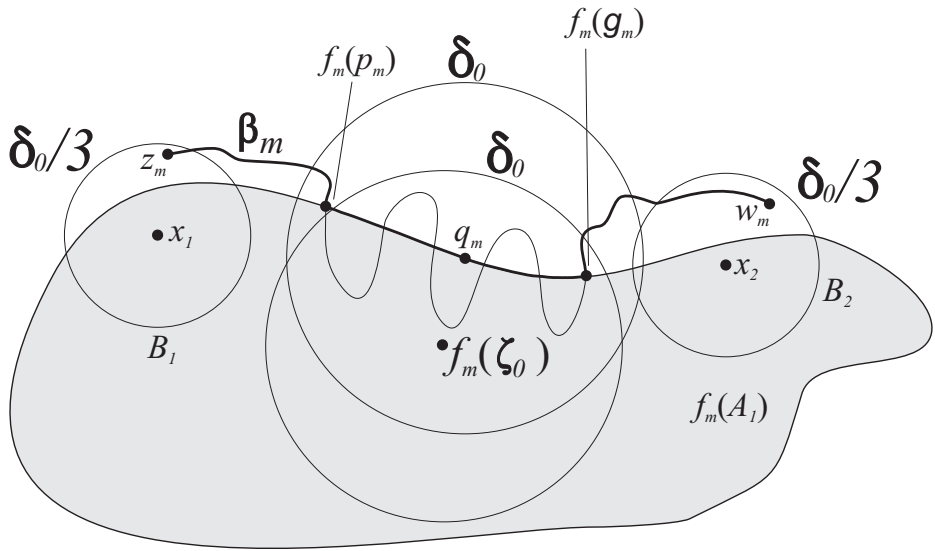}
  \caption{To the proof of Lemma~\ref{lem2}}\label{fig1}
 \end{figure}
We show that points $z_m, w_m$ may be joined by a path $\beta_m$
such that $\beta_m(t)\not\in f_m (A_1)$ for all $t\in [0,1] $ and,
moreover, $\beta_m(t^m_0)\in \partial f_m(A_1)\cap B(f_m(\zeta_0),
\delta_0)$ for some $t^m_0\in (0, 1),$ where $B(f_m(\zeta_0),
\delta_0)=\{x\in {\Bbb M}^n_*: d_*(x, f_m(\zeta_0))<\delta_0\}.$ By
virtue of the Ahlfors regularity of the space ${\Bbb M}_*^n,$ the
condition
\begin{equation}\label{eq13}
v_*(B(f_m(\zeta_0), \delta_0))\geqslant  \widetilde{C}\cdot
\delta_0^n\,,
\end{equation}
is fulfilled with some constant $\widetilde{C}>0$ in (\ref{eq10}).
Due to (\ref{eq10}) and (\ref{eq13}) there is at least one point
$p_m\in ({\Bbb M}_*^n\setminus f_m(A_1))\cap B(f_m(\zeta_0),
\delta_0).$ Since the manifold ${\Bbb M}_*^n$ is connected, the
points $f_m(\zeta_0)$ and $p_m$ can be joined by the path
$\Delta_m:[0, 1]\rightarrow {\Bbb M}_*^n$ such that
$\Delta_m(0)=f_m(\zeta_0),$ $\Delta_m(1)=p_m;$ in this case, since
the space ${\Bbb M}_*^n$ is geodesic (see Lemma~\ref{lem3}), the
path $\Delta_m $ may be chosen so that its length is
$d_*(f_m(\zeta_0), p_m),$ where $d_*(f_m(\zeta_0), p_m)<\delta_0.$
Therefore, $|\Delta_m|\subset B(f_m(\zeta_0), \delta_0).$ Since the
path $\Delta_m$ is not entirely in ${\Bbb M}_*^n\setminus f_m(A_1),$
or $f_m(A_1),$ by~\cite[Theorem~1.I.5.46]{Ku} there is a point
\begin{equation}\label{eq15}
q_m\in
\partial f_m(A_1)\cap |\Delta_m|\subset B(f_m(\zeta_0)\,.
\delta_0)
\end{equation}
Now join the points $z_m$ and $w_m$ by the path $\gamma_m:[0,
1]\rightarrow {\Bbb M}^n_*:$ $\gamma_m(0)=z_m,$ $\gamma_m(1)=w_m.$
At the same time, the path $\gamma_m $ may be chosen so that
$\gamma_m(1/2)=q_m$ due to the connectedness of the manifold ${\Bbb
M}^n_*.$ If, in this case, the path $\gamma_m$ does not belong to
$f_m(A_1),$ we may put $\beta_m:=\gamma_m.$

\medskip
However, let there be at least one point $p\in |\gamma_m|\cap
f_m(A_1).$ Then by~\cite[Theorem~1.I.5.46]{Ku} there is $\omega\in
|\gamma_m|\cap \partial f_m(A_1).$ We put
$$a_m=\inf\limits_{t\in [0, 1], \gamma_m(t)\in \partial f_m(A_1)}t\,,
\qquad b_m=\sup\limits_{t\in [0, 1], \gamma_m(t)\in \partial
f_m(A_1)}t\,.$$
Obviously, $0<a_m\leqslant b_m<1$ and
$\gamma_m(a_m),\gamma_m(b_m)\in \partial f_m(A_1).$ Since $f_m$ is a
homeomorphism of a domain $G$ onto some domain $f_m(G)$ and
$\overline{B(\zeta_0, \varepsilon_0)}\subset G,$ there are unique
points $p_m,$ $k_m$ and $g_m\in S(\zeta_0, \varepsilon_1)$ such that
$f_m(p_m)=\gamma_m(a_m),$ $f_m(k_m)=q_m$ and
$f_m(g_m)=\gamma_m(b_m).$ By virtue of the remarks made at the
beginning of the proof of the lemma, the sphere  $S(\zeta_0,
\varepsilon_1)$ is a path connected set, therefore the points $p_m,$
$k_m$ and $g_m$ may be pairwise joined by a path $\alpha_m:[a_m,
b_m]\rightarrow {\Bbb M}^n$ such that $\alpha_m(t)\in S(\zeta_0,
\varepsilon_0),$ $t\in [a_m, b_m],$ while
$\alpha_m(t_0^m)=k_m=f_m^{\,-1}(q_m)$ for some $t_0^m\in[a_m, b_m].$
Therefore, the path
$$\beta_m(t)= \left\{
\begin{array}{rr}
\gamma_m(t), & t\in [0, a_m]\cup[b_m, 1]\,,\\
f_m(\alpha_m(t)),  &  t\in [a_m, b_m]
\end{array}
\right. $$
is the desired path, namely, $\beta_m$ join the points $w_m$ and
$z_m$ in ${\Bbb M}_*^n\setminus f_m(A_1),$ besides, $\beta_m$ takes
$q_m$ from~(\ref{eq15}) for some $t=t^m_0.$

\medskip
\textbf{V.} So, the path $\beta_m(t)$ is constructed. Now, based on
the path $\beta_m(t),$ we construct some new path
$\overline{\beta}_m$ as follows. If the set $|\beta_m|$ entirely
belongs to $B(q_m, \delta_0),$ we put $\overline{\beta}_m:=\beta_m.$
Note that $d_*(|\beta_m|)\geqslant \delta_0/3,$ which follows
from~(\ref{eq11}). Suppose that there is a point $t^m_1\in [0, 1]$
such that $\beta_m(t^m_1)\not\in B(q_m, \delta_0).$ Then, by the
connectedness of the path $\beta_m,$ there is a point $t^{*m}_2\in
[0, 1]$ such that $\beta_m(t^{*m}_2)\in S(q_m, \delta_0),$
see~\cite[Theorem~1.I.5.46]{Ku}. Then we define $\overline{\beta}_m$
as a subpath of $\beta_m,$ joining the points $q_m$ and
$\beta_m(t^{*m}_2)$ inside the ball $\overline{B(q_m, \delta_0)}.$
Obviously, the geodesic diameter $d _*$ of $\overline{\beta}_m$ is
not less than $\delta_0.$

\medskip
Using reparameterization, if necessary, we may assume that the path
$\overline{\beta}_m$ is defined for $t\in[0, 1].$ Thus, we have
established the existence of a path
$\overline{\beta}_m:[0,1]\rightarrow {\Bbb M}_*^n,$ which has the
following properties:

\medskip
1) the condition $\overline{\beta}_m(t)\in {\Bbb M}^n_*\setminus
f_m(A_1)$ holds for all $t\in [0, 1];$

\medskip
2) the following lower estimate holds:
\begin{equation}\label{eq20}
d_*(|\overline{\beta}_m|)\geqslant \delta_0/3\,,
\end{equation}
where $d_*(|\overline{\beta}_m|),$ as usual, means the geodesic
diameter of the set $|\overline{\beta}_m|\subset {\Bbb M}^n_*;$

\medskip
3) there is an inclusion:
\begin{equation}\label{eq21}
|\overline{\beta}_m|\subset  \overline{B(q_m, \delta_0)}\,.
\end{equation}
Note that the inclusion~(\ref{eq21}) implies another more important
relation, namely, let $x\in |\overline{\beta}_m|,$then by the
triangle inequality $d_*(x, f_m(\zeta_0))\leqslant d_*(x,
q_m)+d_*(q_m, f_m(\zeta_0))\leqslant 2\delta_0,$ see~(\ref{eq15}).
Thus,
\begin{equation}\label{eq22}
|\overline{\beta}_m|\subset  \overline{B(f_m(\zeta_0),
2\delta_0)}\,.
\end{equation}

\medskip
\textbf{VI.} Further reasoning is related to the use of
Proposition~\ref{pr2}. From the lower estimate of the $p$-modulus,
we obtain a contradiction with the assumption made in~(\ref{eq4A}).
Denote $\varepsilon_m:=d(\zeta_m, \zeta_0).$ Recall that
$\varepsilon_m\rightarrow 0$ as $m\rightarrow\infty.$  Join the
points $\zeta_0$ and $\zeta_m$ involved in (\ref{eq4A}) by the path
$\kappa_m:[0, 1]\rightarrow \overline{B(\zeta_0, \varepsilon_m)},$
which is possible because closed balls $\overline{B(\zeta_0,
\varepsilon_m)}$ belong to some normal neighborhood of $\zeta_0,$
and, therefore, $\overline{B(\zeta_0, \varepsilon_m)}$ are path
connected for all $m\in{\Bbb N}.$ Let
$$t_m:=\{\sup t: t\in [0, 1], f_m(\kappa_m(t))\in
B(\zeta_0, \delta_0)\}\,.$$
Let $\eta_m=f_m(\kappa_m)|_{[0, t_m]}$ be part of a path
$f_m(\kappa_m),$ located in the ball $\overline{B(f_m(\zeta_0),
\delta_0)}.$ Note that $|\eta_m|$ is a continuum in
$\overline{B(f_m(\zeta_0), \delta_0)}\subset{\Bbb M}^n_*,$ and that
\begin{equation}\label{eq17}
d_*(|\eta_m|)\geqslant \delta_0\,,
\end{equation}
where, as usual, $d_*(|\eta_m|)$ denotes the geodesic diameter of
the set $|\eta_m|$ in ${\Bbb M}^n_*,$ and $|\eta_m|$ is the locus
(image) of the path $\eta_m.$

\medskip
Now we apply Proposition~\ref{pr2} for $E=|\eta_m|,$
$F=\overline{\beta}_m,$ $X={\Bbb M}^n_*$ and $R=2\delta_0.$ By
(\ref{eq20}), (\ref{eq22}) and (\ref{eq17}), we obtain that $F,
E\subset \overline{B(f_m(\zeta_0, 2\delta_0))},$ $d_*(E)\geqslant
\delta_0$ и $d_*(F)\geqslant \delta_0/3.$ Therefore, by
Proposition~~\ref{pr2}
\begin{equation}\label{eq23}
M_p(\Gamma(E, F, {\Bbb M}^n_*))\geqslant
\frac{1}{\widetilde{C}}\cdot\frac{\min\{d_*(E),
d_*(F)\}}{(2\delta_0)^{1+p-n}}\geqslant
\frac{1}{\widetilde{C}}\cdot\frac{\min\{\delta_0,
\delta_0/3\}}{(2\delta_0)^{1+p-n}}:=C_0>0\,.
\end{equation}
On the other hand, by~\cite[Theorem~1.I.5.46]{Ku} and using the
properties of the modulus of the families of paths, we obtain that
$$\Gamma(F, f_m(\overline{B(\zeta_0, \varepsilon_m)}), {{\Bbb M}^n_*})>
\Gamma(\partial f_m(A_1), f_m(\overline{B(\zeta_0, \varepsilon_m)}),
f_m(A_1))\,,$$
because $F=\overline{\beta}_m(t)\in {\Bbb M}^n_*\setminus f_m(A_1)$
for $t\in [0, 1].$ In addition, since by construction
$E=|\eta_m|\subset f_m(\overline{B(\zeta_0, \varepsilon_m)}),$ then
$$\Gamma(E, F, {\Bbb M}^n_*)\subset\Gamma(F, f_m(\overline{B(\zeta_0, \varepsilon_m)}), {{\Bbb M}^n_*})\,.$$
We obtain from the last two relations and from~(\ref{eq23}) that
\begin{equation*}\label{eq24}
0<C_0\leqslant M_p(\Gamma(E, F, {{\Bbb M}^n_*}))\leqslant
M_p(\Gamma(\partial f_m(A_1), f_m(\overline{B(\zeta_0,
\varepsilon_m)}), f_m(A_1)))=
\end{equation*}
$$={\rm cap}_p\,(f(B(\zeta_0,
\varepsilon_1)), f_m(\overline{B(\zeta_0, \varepsilon_m)}))\,,$$
where the constant $C_0$ depends only on $\delta_0, \widetilde{C},
n$ and $p.$ However, the last relation contradicts~(\ref{eq16A}).
The resulting contradiction indicates the incorrectness of the
assumption in~(\ref{eq4A}), which proves the lemma.~ $\Box$
\end{proof}

\medskip
On the basis of Lemma~\ref{lem2}, the main results of the section
can be formulated.  Their connection with Lemma~\ref{lem2} is
established on the basis of the approach outlined in the proof
of~\cite[Theorems~1.1, 2.1 and 2.2]{IS$_2$},
cf.~\cite[Proposition~2]{Sev$_3$}, and therefore details are
omitted.

\medskip
\begin{theorem}\label{th2}{\sl\,Let $Q:G\rightarrow (0, \infty]$ be a function measurable with
respect to measure $v,$ $n-1<p<n,$ and the Riemannian manifold
${\Bbb M}_*^n$ be considered as a metric space with a geodesic
metric $d_*,$ is $(n, p)$-admissible target. In addition, assume
that ${\Bbb M}_*^n$ is a proper space in which the isoperimetric
inequality~(\ref{eq3}) is fulfilled. Suppose that $Q\in
FMO(\zeta_0).$ Then the family $\frak{F}^p_{\zeta_0, Q}(G, {\Bbb
M}^n_*)$ is equicontinuous at $\zeta_0,$ where equicontinuity is
understood in the sense of the space $({\Bbb M}^n_*, d_*).$
}
 \end{theorem}

\medskip
\begin{theorem}\label{th3}{\sl\,
The conclusion of Theorem~\ref{th2} remains valid if under the
conditions of this theorem, instead of the requirement $Q\in
FMO(\zeta_0),$ we assume that the relation
$\int\limits_{\varepsilon}^{\delta(x_0)}\
\frac{dr}{r^{\frac{n-1}{p-1}}q_{x_0}^{\frac{1}{p-1}}(r)}<\infty$
holds for some $\delta(x_0)>0$ and sufficiently small
$\varepsilon>0,$ while
\begin{equation}\label{eq4}
\int\limits_{0}^{\delta(x_0)}\
\frac{dr}{r^{\frac{n-1}{p-1}}q_{x_0}^{\frac{1}{p-1}}(r)}=\infty
\end{equation}
and $q_{x_0}(r):=\frac{1}{r^{n-1}}\int\limits_{S(x_0,
r)}Q(x)\,d\mathcal{A}.$
}
\end{theorem}

 \medskip
\begin{theorem}\label{th4}{\sl\,
The conclusion of Theorem~\ref{th2} remains valid if under the
conditions of this theorem, instead of the requirement $Q\in
FMO(\zeta_0),$ we assume that $Q\in L_{\rm loc}^s({\Bbb R}^n),$
where $s\geqslant\frac{n}{n-p}.$
}
 \end{theorem}

\section{Examples}
\setcounter{equation}{0}

\begin{example}\label{ex1}
Consider, first of all, the family of analytic functions
$f_n(z)=e^{nz},$ $n\in {\Bbb N}.$ Note that this family of mappings
is not equicontinuous at the point $z_0=0.$ In this regard, let us
specify which of the conditions of Theorems~\ref{th1} and~\ref{th2}
are violated here. First of all, $f_n$ are ring $Q_n$-maps in the
whole complex plane for $Q_n(z)\equiv 1,$ see
e.g.~~\cite[Theorem~2]{Pol}, cf.~\cite[Theorems~8.1, 8.6]{MRSY}.
Moreover, $f_n$ are a ring $Q_n$-maps with respect to $p$-modulus,
$1<p\leqslant 2,$ for $K_{I, p}(z,
f_n)=Q_n(z)=n^{2-p}|e^{nz}|^{2-p},$ see~\cite[Theorem~1.1]{SalSev}.
Here, of course, the role of measures $\mu$ and $\mu^{\,\prime}$ is
performed by the Lebesgue measure on the complex plane. The function
$K_{I, p}(x, f)$ used above, is calculated by the formula
$$K_{I, p}(x,f)\quad =\quad\left\{
\begin{array}{rr}
\frac{|J(x,f)|}{{(l(f^{\,\prime}(x)))^p}}, & J(x,f)\ne 0,\\
1,  &  f^{\,\prime}(x)=0, \\
\infty, & {\rm otherwise}
\end{array}
\right.\,,$$
where $J(x,f)=\det f^{\,\prime}(x)$ and
$l(f^{\,\prime}(x))=\min\limits_{h\in {\Bbb R}^n \backslash \{0\}}
\frac {|f^{\,\prime}(x)h|}{|h|}.$

Note that, in the case $p=2,$ the function $Q\equiv 1$ has a finite
mean oscillation at each point of the complex plane. However,
Theorem~\ref{th1} may not be applicable to the family of mappings
$\{f_n\}_{n=1}^{\infty},$ since these maps do not omit any fixed
continuum $K$ in ${\Bbb C}.$ The situation does not change, if we
consider the same family of mappings in a certain bounded subdomain
of the complex plane.

Since the omitting of some continuum $K$ is not required under the
conditions of Theorem~\ref{th2}, we consider separately the case
$p<1<2.$ One could put here ${\Bbb M}^2={\Bbb M}^2_*={\Bbb C}$ and
$d(x, y)=d_*(x, y)=|x-y|,$ $x, y\in {\Bbb C}.$ However, in this
case, another condition of this theorem is violated, namely, the
functions $Q_n$ do not have a common majorant $Q$ that has a finite
mean oscillation at $0.$ It is easy to see that all the other
conditions of Theorem~\ref{th2} are satisfied.
\end{example}

\medskip
\begin{example}\label{ex2}
For fixed $1<q<2,$ consider the family of mappings
$f_m(x)=\left(\frac{x_1}{m}, \frac{x_2}{m^q}\right),$ $m=1,2,\ldots
\,,$ $z=(x_1, x_2)\in {\Bbb D},$ ${\Bbb D}=\{z\in {\Bbb R}^2:
|z|<1\}.$ Notice, that $f_m$ are homeomorphisms of the unit disk
into itself, wherein, $J(z, f_m)=\frac{1}{m^q}$ and
$l(f_m^{\,\prime}(x))=\frac{1}{m}.$ By~\cite[Theorem~1.1]{SalSev}
$f_m$ a ring $Q_m$-homeomorphisms with respect to $p$-modulus for
$Q_m=K_{I, p}(z, f_m)=m^{p-q}.$ It is clear from the definition that
the family $\{f_m\}_{m=1}^{\infty}$ is equicontinuous, say, at the
point $0,$ however, this conclusion cannot be obtained from
Theorem~\ref{th1} or~\ref{th2} if we apply these theorems for the
exponent $p=\alpha=\alpha^{\,\prime}=2.$ Namely, the functions $Q_m$
obviously do not have a common majorant of the corresponding class.

However, this result can be obtained from the same theorems by
setting $p$ to be less than $q.$ Then we may put $Q(x)\equiv 1,$
since $Q_m(x)\leqslant 1$ everywhere in ${\Bbb D}.$
\end{example}

\medskip
\begin{example}\label{ex3}
Consider the family of mappings from Example~\ref{ex2}, with the
difference that we define mappings $f_m$ in the entire space ${\Bbb
R}^2.$ In this case, the equicontinuity of the family
$\{f_m\}_{m=1}^{\infty}$ at 0 cannot be directly obtained from
Theorem~\ref{th1}, since the mappings do not omit any continuum $K$
in ${\Bbb R}^2.$ However, this result directly follows from
Theorem~\ref{th2}, in this case, since the maps $f_m$ are
homeomorphisms.
\end{example}

\medskip
\begin{example}\label{ex5}
\medskip
Fix a number $p\geqslant 1$ satisfying the condition $n/p(n-1)<1.$
Put $\alpha\in (0, n/p(n-1)).$ We define the sequence of
homeomorphisms~$f_m$ of the unit ball onto $B(0, 2)$ as follows:
$$f_m(x)\,=\,\left
\{\begin{array}{rr} \frac{1+|x|^{\alpha}}{|x|}\cdot x\,, & 1/m\leqslant|x|< 1, \\
\frac{1+(1/m)^{\alpha}}{(1/m)}\cdot x\,, & 0<|x|< 1/m \ .
\end{array}\right.
$$
Put ${\Bbb B}^n=\{x\in {\Bbb R}^n: |x|<1\}.$ Note that $f_m$ a ring
$Q$-homemorphisms with respect to $n$-modulus for
$Q=\left(\frac{1+|x|^{\,\alpha}}{\alpha
|x|^{\,\alpha}}\right)^{n-1}$ and for each $\zeta_0\in
\overline{{\Bbb B}^n},$ moreover, $Q\in L^p({\Bbb B}^n),$ see,
e.g.~\cite[proof of Theorem~7.1]{Sev$_4$}. Note that $\frak
F=\{f_m\}_{m=1}^{\infty}$ is not equicontinuous in ${\Bbb B}^n.$
Indeed, $|f_m(x_m)-f(0)|=1+1/m \rightarrow 1$ as
$m\rightarrow\infty,$ where $|x_m|=1/m.$ The reason for the latter
circumstance, as well as the impossibility of applying
Theorems~\ref{th1} and~\ref{th2}, is that $Q$ does not belong to the
$FMO$ class at $\zeta_0=0.$
\end{example}

\medskip
\begin{example}\label{ex4}
Let ${\Bbb D}$ be the unit disk on the plane. The hyperbolic
distance in ${\Bbb D}$ is given by the formula
\begin{equation*}\widetilde{h}(z_1, z_2)=\log\,\frac{1+t}{1-t}\,,\quad
t=\frac{|z_1-z_2|}{|1-z_1\overline{z_2}|}\,,
\end{equation*}
while the hyperbolic area of a set $S$ in ${\Bbb D}$ is calculated
as the integral
$\widetilde{h}(S)=\int\limits_S\frac{4\,dm(z)}{(1-|z|^2)^2},$ see
e.g.~\cite[(2.4)--(2.5)]{RV}.
Given a Borel function $\rho:{\Bbb D}\rightarrow [0, \infty],$ a
Lebesgue measurable set $S\subset {\Bbb D}$ and a locally
rectifiable path $\gamma:(a, b)\rightarrow {\Bbb D}$ we set
\begin{equation*}
\int\limits_S\rho(z)\,d\widetilde{h}(z):=\int\limits_S\frac{4\rho(z)\,dm(z)}{(1-|z|^2)^2}\,,\quad
\int\limits_{\gamma}\rho(z)\,ds_{\widetilde{h}}(z):=\int\limits_{\gamma}\frac{2\rho(z)\,|dz|}{1-|z|^2}\,.
\end{equation*}
Note that the metric space space $({\Bbb D}, \widetilde{h})$ is not
regular by Ahlfors, although it is regular from below,
see~\cite[Example~8.24(c) and Proposition~8.19]{He}. Fix $r_0>0$ and
note that the space $(B(0, r_0), \widetilde{h})$ with hyperbolic
area is Ahlfors 2-regular, where $B(z_0, r)=\{z\in {\Bbb C}:
|z-z_0|<r\}.$ Indeed, using the notation $B_{\widetilde{h}}(z_0,
r)=\{z\in {\Bbb D}: \widetilde{h}(z, z_0)<r\},$ we obtain that
\begin{equation}\label{eq1A}
B(z_0, C_2r)\subset B_{\widetilde{h}}(z_0, r)\subset B(z_0, C_1r)\,,
\end{equation}
where $C_1$ and $C_2$ are some positive constants depending on
$r_0.$ We also note that the hyperbolic area $\widetilde{h}$ is
related to the Lebesgue measure on the complex plane by the
relations
\begin{equation}\label{eq3D}
C_4\cdot m(E)\leqslant\widetilde{h}(E)\leqslant C_3\cdot m(E)
\end{equation}
whenever $E\subset B(z_0, r_0)$ and $C_3$ and $C_4$ are some
positive constants depending on $r_0.$ In this case, the Ahlfors
regularity of the space $(B(0, r_0), \widetilde{h})$ follows from
relations~(\ref{eq1A}) and (\ref{eq3D}), as well as the fact that
ordinary Euclidean balls are Ahlfors regular,
see~\cite[Example~8.24(a) and Proposition~8.19]{He}.

\medskip
Now consider the following family of mappings
$$f_m(x)\,=\,\left
\{\begin{array}{rr} \frac{1+|x|^{\alpha}}{|x|}\cdot x\,, & 1/m\leqslant|x|<r_0, \\
\frac{1+(1/m)^{\alpha}}{(1/m)}\cdot x\,, & 0<|x|< 1/m \ .
\end{array}\right.
$$
We will consider $f_m,$ $m=1,2,\ldots ,$ as mappings acting between
spaces $(B(0, r_0), \widetilde{h})$ and $({\Bbb C}, |\cdot|),$ where
$|\cdot|$ corresponds to the Euclidean distance, $(B(0, r_0),
\widetilde{h})$ equipped with the hyperbolic area, and $({\Bbb C},
|\cdot|)$ equipped with the Lebesgue measure. According
to~\cite[Remark~5.2]{Sev$_5$}, the mappings $f_m$ are ring $Q$-maps
for $Q=\left(\frac{1+|x|^{\,\alpha}}{\alpha
|x|^{\,\alpha}}\right)^{n-1}$ at $0$ with respect to 2-modulus.
Here we are talking about ring $Q$-mappings in terms of the
specified spaces $(B(0, r_0), \widetilde{h})$ and $({\Bbb C},
|\cdot|),$ that is, with respect to the corresponding metrics and
measures. As we noted in Example~\ref{ex5}, the family of mappings
$f_m$ is not equicontinuous. In the same time, all conditions on the
spaces involved in Theorem~\ref{th1}, are fulfilled, including the
Ahlfors regularity of the space~$(B(0, r_0), \widetilde{h}).$ It can
be shown that $Q\not\in FMO(0)$ in terms of the metric and measure
of this space.
\end{example}

\medskip
\begin{example}\label{ex6}
To obtain a family of mappings that are equicontinuous at $0,$ in
Example~\ref{ex5} one could set
$$
g_m(x)=\left\{
\begin{array}{rr}
\frac{x}{|(m-1)/m|\log\frac{e}{(m-1)/m}}, & x\in {\Bbb B}^n\cap B(0, (m-1)/m),\\
\frac{x}{|x|\log\frac{e}{|x|}}, &  x\in {\Bbb B}^n\setminus B(0,
(m-1)/m)
\end{array}
\right.\,,
$$
where ${\Bbb B}^n=\{x\in {\Bbb R}^n: |x|<1\}.$ It can be shown that
$g_m$ are ring $Q$-maps with respect to the $n$-module, where
$Q(x)=\log^{n-1}\frac{e}{|x|},$ see the reasoning used in the
consideration of~\cite[Proposition~6.3]{MRSY}. The equicontinuity of
this family at zero can be obtained by direct calculations, however,
it also follows from~Lemma~\ref{lem2}. In particular, the existence
of the corresponding function $\psi$ follows
from~\cite[Lemma~7.3]{MRSY}, while the divergence of the integral of
the form~(\ref{eq4}) for $Q(x)=\log^{n-1}\frac{e}{|x|}$ is directly
verified.

It is also easy to indicate a similar family of mappings inside the
unit disk with a hyperbolic metric and a hyperbolic area.
\end{example}

%=================Список литературы====================
%\end{fulltext}

\medskip
\medskip
{\bf \noindent Evgeny Sevost'yanov, Sergei Skvortsov} \\
Zhytomyr Ivan Franko State University,  \\
40 Bol'shaya Berdichevskaya Str., 10 008  Zhytomyr, UKRAINE \\
Email: esevostyanov2009@gmail.com, serezha.skv@gmail.com

\medskip
\noindent{{\bf Evgeniy Petrov}\\
Institute of Applied Mathematics and Mechanics\\
of NAS of Ukraine, \\
1 Dobrovol'skogo Str., 84 100 Slavyansk,  UKRAINE\\
Email: eugeniy.petrov@gmail.com

\end{document}